\documentstyle[art12,amsfonts]{article}
\parskip2pt\renewcommand{\sl}{\it}\newcommand{\ns}{\normalshape}
\newcommand{\n}{\noindent}\newcommand{\rf}[1]{(\ref{#1})}
\mathsurround1pt\newcommand{\k}{\kern}
\newcommand{\ds}{\displaystyle}\newcommand{\ts}{\textstyle}
\newcommand{\lra}{\longrightarrow}
\newcommand{\chii}{\raise2pt\hbox{$\chi$}}\newcommand{\W}{\Omega}
\newcommand{\hpq}{h^{p,q}} \renewcommand{\sp}[1]{{}\!^{#1}}
\newcommand{\vep}{\varepsilon}

\newcommand{\Ga}{\Gamma}\newcommand{\p}{\hbox{\ns\bf p}}
\renewcommand{\c}{\hbox{\ns\bf c}}\newcommand{\s}{\hbox{\ns\bf s}}
\newcommand{\K}{\hbox{\ns\bf K}} \newcommand{\M}{{\cal M}}
\newcommand{\hol}{\hbox{\ns Hol}}\newcommand{\w}{\omega}
\def\R{{\Bbb R}\k.5pt}\newcommand{\C}{{\Bbb C}\k.5pt}
\newcommand{\Sp}{\hbox{\ns Sp}}\newcommand{\CP}{{\Bbb C\Bbb P}}
\newcommand{\Z}{{\Bbb Z}}\newcommand{\spin}{\hbox{\ns Spin}(7)}
\newcommand{\G}{\hbox{\ns G}}
\newcommand{\SO}{\hbox{\ns SO}}\newcommand{\SU}{\hbox{\ns SU}}
\newcommand{\hk}{hyper-K\"ahler}\newcommand{\ka}{K\"ahler}
\newcommand{\ch}{\hbox{\ns\bf ch}}\newcommand{\td}{\hbox{\ns\bf td}}
\renewcommand{\a}{&\k-3pt=\k-3pt&\ds}\newcommand{\D}{\partial}
\newcommand{\ga}{\gamma}\newcommand{\bin}[2]{\hbox{${#1\choose#2}$}}
\newcommand{\ol}{\overline}\newcommand{\ul}{\underline}

\newcommand{\ext}[1]{\raise1pt\hbox{$\ts\bigwedge$}\k-1pt^{#1}\k-1pt}
\newcommand{\ba}{\begin{array}}
\renewcommand{\bar}{\begin{array}{rcl}}
\newcommand{\be}{\begin{equation}}
\newcommand{\ea}{\end{array}}\newcommand{\ee}[1]{\label{#1}
\end{equation}}\newcommand{\frt}[2]{\hbox{$\ts\frac{#1}{#2}$}}
\newcommand{\frs}[2]{\hbox{\large$\ts\frac{#1}{#2}$\normalsize}}\outer
\def\pro#1#2\par{\bigbreak\noindent{\bf#1.\enspace}{\sl#2}\par\bigbreak}
0 
\begin{document}

\centerline{\LARGE\bf Cohomology of K\"ahler Manifolds}
\vspace{7pt}\par\centerline{\LARGE\bf with $\c\lower2.5pt
\hbox{$_{\ns\bf1}$}$=\kern2.5pt0}\par\vspace{.25in}
\centerline{\Large S.$\,$M.\ Salamon}\par\vspace{.25in}
\centerline{\it To Professor Calabi on his 70th birthday}
\vspace{.3in}

\subsection*{Introduction}

A more informal title of this report might be `Betti, Hodge and Chern
numbers', without an implied order of preference. To begin with, let
$M$ be a compact connected \ka\ surface, so that the real dimension of
$M$ is 4. The Betti numbers of $M$ are given in terms of Hodge numbers
by the usual relations \[ \ba{c}b_1=h^{1,0}+ h^{0,1}= 2h^{0,1} \\[3pt]
b_2=h^{2,0}+h^{1,1}+ h^{0,2}=h^{1,1}+2h^{0,2}.\ea\] Noether's formula
states that
\[\big<\frt1{12}(\c_1\sp2+\c_2),[M]\big>=1-h^{0,1}+h^{0,2},\] and the
left-hand side is by definition the Todd genus of $M$. Moreover, the
Euler characteristic $e(M)$ is equal to either side of the equation \[
\left<\c_2,[M]\right>=2-2b_1+b_2.\] \smallbreak

If the first Chern class $\c_1$ of $M$ vanishes over $\Bbb Z$, then
the canonical bundle of $M$ is trivial, and $h^{0,2}=1$. As a
consequence, we obtain the formula \be 0=22 - 4b_1 - b_2.\ee{4} In
these circumstances, the classification of complex surfaces implies
that $M$ must in fact be a torus $T$ ($b_1=4$, $b_2=6$) or a K3
surface $K$ ($b_1=0$, $b_2=22$). We may construct $K$, at least up to
diffeomorphism, by resolving the 16 singular points of the orbifold
$T/\Z_2$, where $\Z_2$ is generated by the mapping
$x\mapsto-x$. Observe that the right-hand side of \rf{4} assumes the
value 16 for $T/\Z_2$ ($b_1=0$, $b_2=6$); in fact, we shall see that
the formula has a significance which transcends the examples.

We move on to consider the case in which $M$ is now a \ka\ manifold of
{\sl complex} dimension 4. The introduction to \cite{H} contains a
formula which (after slight rearrangement) bears a striking similarity
to \rf{4}. Namely, if $\chii^p$ denotes the alternating sum
$\sum_{q=0}^4(-1)^q\hpq$ of Hodge numbers of $M$ then
\be\left<\c_1\c_3,[M]\right>=22\chii^0-4\chii^1-\chii^2.\ee{hir} This
equation is generalized by Theorem~2, which asserts that on any almost
complex $n$-dimensional manifold the characteristic number
$\left<\c_1\c_{n-1},[M]\right>$ can be expressed in terms of the
integers $\chii^p$, $0\le p\le[n/2]$. This is a simple corollary of
the Riemann-Roch theorem that was derived in \cite{coho} (an expanded
version of this article) by differentiating a suitable K-theoretical
expression. The author subsequently discovered that the corollary had
been used earlier by Narasimhan and Ramanan
\cite{NR}, in a context discussed briefly at the end of Section~3. In
any case, we deduce below that if $n=2m$ is even and $\c_1=0$ then the
`mirror symmetry' \be\hpq\longleftrightarrow h^{n-p,q},\quad 0\le
p,q\le n,\ee{ms} reverses the sign of a certain linear combination
$\Phi$ of Betti numbers, defined in Theorem~1. This complements the
obvious fact that the operation \rf{ms} acts as $(-1)^n$ on the Euler
characteristic $e$, and adds weight to the view that \rf{ms} has some
significance in higher dimensions.

Compact \hk\ manifolds constitute a class of \ka\ manifolds with
$\c_1\!=\!0$ and whose Hodge numbers are invariant by the symmetry
\rf{ms}. It follows from above that their Betti numbers are subject to
the equation $\Phi=0$, of which \rf{4} is a special case.  This
provides an effective means to prove the non-existence of \hk\ metrics
on particular manifolds. The name `\hk' was introduced by Calabi
\cite{C2,Cal}, who gave the first non-trivial examples of these
metrics, and highlighted their importance prior to their appearance on
symplectic quotients and on moduli spaces of anti-self-dual
connections \cite{Hi,Kob,Kr,Mu}. Although our `\hk\ constraint'
$\Phi=0$ is ultimately derived from the Atiyah-Singer Index Theorem,
it is in some ways an independent result and may be expected to have a
more direct proof. Quite why it respects functorial properties of \hk\
manifolds is briefly explained in Section~3, with reference to the now
standard examples of Beauville \cite{Bea}; more details are in
\cite{coho}. The remainder of Section~3 is of a different character,
as we discuss the topology of a selection of topical examples which
fall outside the scope of the title, but which nonetheless are
suggestive of a more extensive theory linking geometrical structures
and cohomology.

\vfil\eject

\subsection*{1. Background and main results}

Throughout, $M$ denotes a compact oriented manifold of real dimension
$d$.  Recall that $M$ is {\sl\ka} if it possesses both a Riemannian
metric $g$ and an orthogonal complex structure $J$ for which the
resulting 2-form \be \w(X,Y)=g(JX,Y),\qquad X,Y\in TM,\ee{w} is
closed. We shall always denote the complex dimension of $M$ by $n$, so
that $d=2n$. We denote the $k$th Chern class of (the holomorphic
tangent bundle of) $M$ by $\c_k$, and we shall generally regard it as
an element of $H^{2k}(M,\R)$. In particular, the first Chern class
$\c_1$ of $M$ is represented by $1/2\pi$ times the Ricci tensor
(converted into a 2-form by $J$ in analogy to \rf{w}). If $\c_1$
vanishes over $\R$, then Yau's proof of the Calabi conjecture implies
that $M$ admits a \ka\ metric with zero Ricci tensor \cite{Y,Be}. This
implies that the holonomy group $\hol(M)$ of the Levi\,Civita
connection is contained in the subgroup $\SU(n)$ of $\SO(d)$.

For example, a quartic hypersurface in $\CP^3$ is a K3 surface which
is manifestly \ka\ and has $\c_1=0$; it must therefore admit a \ka\
(`Calabi-Yau') metric with holonomy group equal to $\SU(2)$. The
latter is best identified with the group $\Sp(1)$ of unit
quaternions. More generally, the Cheeger-Gromoll splitting theorem
implies that a compact \ka\ manifold with $\c_1=0$ has a finite
holomorphic covering of the form \be T^{2k}\times
X_1\times\cdots\times X_r\times Y_1\times \cdots\times Y_s,\ee{cov}
where $T^{2k}\cong\C^k/\Z^{2k}$ is a complex torus, $X_i$ and $Y_j$
are simply-connected, and \[\hol(X_i)=\SU(n_i),\qquad
\hol(Y_j)=\Sp(m_j)\] with $m_j\ge1$ and $n_i\ge3$.  Further details
may be found in \cite{Bea,Be,Bog,LM}.

A Riemannian manifold $M$ of real dimension $d=4m$ which, like the
$Y_j$ or their products, has $\hol(M)\subseteq\Sp(m)$ is said to be
{\sl hyper-K\"ahler}. Such a manifold can be characterized by the
existence of a a triple of orthogonal almost complex structures
$J_1,J_2,J_3$ satisfying the algebraic condition $J_1J_2=J_3=-J_2J_1$
and such that the associated 2-forms $\w_1,\w_2,\w_3$ defined by
\rf{w} are all closed. The latter condition (in contrast to the case
of a single almost comnplex structure) implies that each $J_i$ gives
rise to a complex, and therefore \ka, structure \cite{Hi,parma}. A
\hk\ manifold $M$ actually possesses a continuous
family of complex structures, namely $\sum_{i=1}^3 a_iJ_i$ with
$\sum_{i=1}^3(a_i)^2=1$, parameterized by $S^2$ or $\CP^1$. The
resulting deformation or twistor space $M\times\CP^1$ was emphasized
from the start in Calabi's study \cite{Cal}. A \hk\ manifold $M$ is
irreducible if and only if $\hol(M)=\Sp(m)$; in this case $h^{2,0}=1$.
By contrast, the manifolds $X_i$ in \rf{cov} have $h^{2,0}=0$.

We shall denote the Poincar\'e polynomial of the compact oriented
$d$-manifold $M$ by $P(M;t)$, or by $P(t)$ if the latter causes no
confusion. Thus, if the Betti numbers of $M$ are denoted by $b_j$,
$0\le j\le d$, then \be P(t)=\sum_{j=0}^d b_j t^j=\sum_{j=0}^d b_j
t^{d-j}.\ee{P} The alternating sum $P(-1)$ equals the Euler
characteristic $e(M)$; our first main result concerns an analogous
expression.

\pro{Theorem 1} Let $M$ be a compact oriented manifold of real
dimension $d=4m$. Set \[\Phi(M)= 6P''(-1)+\frt12d(5-3d)P(-1).\] If $M$
has a \hk\ metric then $\Phi(M)=0$.

\n If $m=1$ then $\Phi(M)$ equals twice the right-hand
side of \rf{4}; for $m=2$ we obtain

\be  \frs14\Phi(M)= 46-25b_1+10b_2-b_3-b_4.\ee{8}

\n A proof of Theorem~1 based upon index theory for the
Dirac operator was sketched in \cite{torino}. Below, we shall derive
it from a corresponding result for almost complex manifolds, Theorem~2
below, by representing $\Phi$ as a suitable characteristic number.

There is an elementary reason why the first derivative $P'(-1)$ is not
needed in the definition of $\Phi(M)$. For from \rf{P}, \[ P'(-1)=
-\sum_{j=1}^d(-1)^jjb_j=-\sum_{j=1}^d(-1)^j(d-j)b_j,\] and adding the
last two members, \be 2P'(-1)=-d\,P(-1).\ee{der} Using this equation,
Theorem~1 may be expressed in the equivalent form \be
m\,e(M)=6\sum_{j=0}^{2m-1}(-1)^j(2m-j)^2b_j.\ee{fo} It is well known
that on a \ka\ manifold, $b_j\equiv0$ mod$\>2$ if $j$ is odd;
similarly on a \hk\ manifold, $b_j\equiv0$ mod$\>4$ if $j$ is odd
\cite{W,F}. The next result combines this fact with \rf{fo}.

\pro{Corollary} Let $M$ be a compact \hk\ $4m$-manifold. Then
$m\,e(M)\equiv0$ {\ns mod}$\>24$, and in particular $e(M)\equiv0$ {\ns
mod}$\>2$ if $m\not\equiv0$ {\ns mod}$\>8$.

Now let $M$ be a compact almost complex manifold of real dimension
$2n$. The choice of an almost Hermitian metric on $M$ enables one to
define the formal adjoint $\ol\D{}^*=-*\kern-1pt\ol\D\kern1pt*$ of the
$\ol\D$ operator. There is then an elliptic differential operator
\[\bigoplus_{q~\ns even}\W^{p,q}\stackrel{\ol\D+\ol\D^*}
{\lra}\bigoplus_{q~\ns odd}\W^{p,q},\] whose index is denoted by
$\chii^p$ in the notation of \cite{H}. The next results are valid in
this general setting, although when the almost complex structure on
$M$ is integrable, $\chii^p$ is more conveniently defined as
$\sum_{q=0}^n(-1)^q\hpq$, where $\hpq$ is the dimension of the
corresponding Dolbeault cohomology space or equivalently of the \v
Cech space $H^q(M,{\cal O}(\ext pT^*))$.

In all cases there is the `Serre duality' relation \be\chii^{n-p}
=(-1)^n\chii^p,\ee{Sd} and we set \be\chii(t)=\sum_{p=0}^n
\chii^pt^p=(-1)^n\sum_{p=0}^n\chii^{n-p}t^p,\ee{set} which is
often denoted by $\chii_t$. For example, \be\chii(-1)=
\sum_{p=0}^n(-1)^p \chii^p =\sum_{j=0}^{2n}(-1)^jb_j\ee{e1} coincides
with $P(-1)=e(M)$. As we shall explain in the next section, the
well-known formula\be\left<\c_n,[M]\right>=\chii(-1)\ee{e2} may be
regarded as the first of a sequence expressing the coefficients of the
polynomial $\chii(-1-t)$ (which is in some ways more natural than
$\chii(t)$) in terms of Chern numbers. The quadratic term yields

\pro{Theorem 2} Let $M$ be a compact almost complex manifold of
real dimension $2n$. Then \[\left<\c_1\c_{n-1},[M]\right> =
6\chii''(-1)+\frt12n(5-3n)
\chii(-1).\] \bigbreak

\n An equivalent version of the formula can be found at the end of
\cite{NR}. If the complex dimension $n$ is odd then \rf{e1} and \rf{e2}
imply immediately that $e(M)$ is even. The following are slightly less
obvious consequences of Theorem~2 implicit in the theory of $\SU$
cobordism described in \cite{St}.

\pro{Corollary} Let $M$ be a compact almost complex manifold of real
dimension $2n$ with $\c_1=0$. Then {\ns(i)} $n\,e(M)\equiv0$ {\ns
mod}$\;3$, and {\ns(ii)} if $n\equiv2$ {\ns mod}$\;4$ then
$e(M)\equiv0$ \hbox{\ns mod}$\;2$.

As regards (i), which is also follows from \cite{H53}, we remark that
when $n=3$ there are some familiar manifolds wit $3\kern-2pt\not
\kern-1.2pt|\kern1.2pt\chii$ admitting non-integrable almost complex
structures, for instance $S^6$ and $\CP^3$. If $M$ is simply-connected
and $\c_1=0$ then the structure group of $M$ reduces to $\SU(n)$ and
lifts to $\hbox{\ns Spin}(2n)$, and (ii) also follows from Ochanine's
theorem \cite{O}. The latter states that a compact oriented smooth
spin manifold of real dimension $d\equiv4$ mod 8 has signature
divisible by 16 (and therefore $e(M)$ divisible by 2). Other
divisibility properties of the Chern numbers appear in the next
section, and related results on Chern (or rather Segre) numbers can
also be found in \cite{Ar} and references therein.

\subsection*{2. Proofs and generalizations}

In this section, $T$ denotes the holomorphic tangent bundle of an
almost complex manifold of real dimension $2n$. We give the total
Chern class of $T$ a formal factorization \[\c(T)=\prod_{i=1}^n
(1+x_i),\] so that the Chern classes of $T$ may be regarded as the
elementary symmetric polynomials in the symbols $x_1,\ldots,x_n$. The
Chern character of $T$ is then given by \be
\ch(T)=\sum_{i=1}^n e^{x_i}=n+\s_1+\frs1{2!}\s_2+\frs1{3!}\s_3
+\cdots\ee{ch} where $\s_k=\sum_{i=1}^nx_i\sp k$. We shall in fact
need $\ch(\ext pT^*)$, which is formed by replacing the $n$ elements
$x_i$ by the $\bin np$ elements $-(x_{i_1}+x_{i_2}+\cdots+x_{i_p})$,
$i_1<i_2<\cdots<i_p$, which are the weights of the representation
defining $\ext pT^*$. Recall also that the Todd class of $T$ is given
by \[\td(T)=\prod_{i=1}^n\frac{x_i}{1-e^{-x_i}} =
1+\frs12\c_1+\frs1{12} (\c_1\sp2+\c_2)+\frs1{24}\c_1\c_3+\cdots\]

With the above preliminaries, the general form of the
Hirzebruch-Riemann-Roch theorem \cite{AS,H} implies, making use of
\rf{set}, that \be\bar\chii(t)\a(-1)^n\sum_{p=0}^n t^p\left<\ch
(\ext{n-p}T^*)\td(T),[M]\right>.\\[16pt]\a (-1)^n\left<\ch(V(t))
\td(T),[M]\right>,\ea\ee{RR} where \be V(t)=\sum_{p=0}^n
t^p\ext{n-p}T^*\ee{poly} is regarded as an element of $K(M)[t]$. Using
the exterior power operation of K-theory, we may write \[ V(-1)=\ext
n(T^*-\ul{\C}),\] where $\ul\C$ denotes a trivial line
bundle. Moreover, \rf{poly} may be differentiated with respect to $t$
to obtain analogous expressions \[ \bar
V'(-1)\a\ext{n-1}(T^*-\ul{\C}^2),\\[4pt]
V''(-1)\a2\ext{n-2}(T^*-\ul\C^3),\ea\] where $\ul\C^k$ denotes a
trivial line bundle with fibre $\C^k$.

Suppose for a moment that $T^*$ contains $\ul\C^3$ as a subbundle, so
that $T^*-\ul\C^3$ is a genuine complex vector bundle of rank
$n-3$. This effectively corresponds to the case in which $\c_k=0$ for
$k>n-3$.  Then $V''(-1)$ is zero, merely by virtue of its
dimension. Accordingly, we may deduce that in general $\ch(V''(-1))$
belongs to the ideal $\left<\c_{n-2},\c_{n-1},\c_n\right>$ generated
by the `top three' Chern classes, and the proof of Theorem~2 is
completed by the more precise

\pro{Lemma} $(-1)^n\ch(V''(-1))\td(T)=2\c_{n-2}+(n\!-\!1)\c_{n-1} +
\frs1{12} (2\c_1\c_{n-1}+n(3n\!-\!5)\c_n)$.

\n This equation is derived from the formal factorization
\be (-1)^n\ch(V(-1-t))\td(T)=\prod_{i=1}^n\Big(x_i+t\frac{x_i}
{1-e^{-x_i}}\Big).\ee{xxx} The coefficient of $t^2$ in \rf{xxx} is the
sum of $\bin n2$ terms, one of which is \be \big[\,1+\frt12(x_1+x_2) +
\frt1{12}(x_1\sp2+3 x_1x_2+x_2\sp2)\,\big] x_3x_4\cdots x_n.\ee{11}
On the other hand, $\c_1\c_{n-1}$ gives rise to a sum of $n$ terms,
one of which is \be (x_1+\cdots+x_n)x_2\cdots x_n = x_1x_2\cdots x_n +
x_2\sp2x_3\cdots x_n+\cdots\ee{22} The proof of the lemma is completed
from a comparison of \rf{11} and \rf{22}. \smallbreak

We may rewrite \rf{RR} in the form \[\chii(-1-t) =
\left<\K_n^\bullet(t), [M]\right>,\] or equivalently \be
\chii^{(k)}(-1)=(-1)^kk!\left<\K_{n,k}^\bullet,[M]\right>,\ee{K}
where \[\K^\bullet_n(t)=\K_{n,0}^\bullet +\K_{n,1}^\bullet
t+\K_{n,2}^\bullet t^2+\cdots\] denotes the component of either side
of \rf{xxx} in $H^{2n}(M,\R)$ (we use the notation of \cite{coho} in
which the bullet indicates a class of top degree). A generalization of
the above lemma asserts that both $\K^\bullet_{n,2k}$ and
$\K^\bullet_{n,2k+1}$ belong to \[H^{2n}(M,\R)\>\cap\>\left<
\c_{n-2k+1},\c_{n-2k+2},\ldots,\c_n\right>,\] with the exception of
$\K_{n,0}^\bullet$ which equals $\c_n$. In fact $\K_{n,2k+1}^\bullet$
is a linear combination of $\K_{n,2j}^\bullet$ for $0\le j\le k$, and
we have the following explicit formulae.

\pro{Lemma}\newcommand{\hs}{\\[3pt]\hspace{45pt}}
{\small\[\ba{l}2^23\,\K_{n,2}^\bullet=\c_1\c_{n-1}+\frt12n(3n-5)\c_n,
\\[10pt]2^43^25\,\K_{n,4}^\bullet=[-\c_1\sp3+3\c_1\c_2-3\c_3]\c_{n-3}+
[\c_1\sp2+3\c_2]\c_{n-2}+\frt12[15n^2-85n+108]\c_1\c_{n-1}\hs
+\frt18n[15n^3-150n^2+485n-502]\c_n\\[10pt]2^53^35^17\,\K_{n,6}^\bullet
=[\c_1\sp5-5\c_1\sp3\c_2+5\c_1\c_2\sp2+5\c_1^2\c_3-5\c_2\c_3-5\c_1\c_4+
5\c_5]\c_{n-5}\hs+\frt12[-2\c_1\sp4+\c_1\sp2\c_2+10\c_2\sp2-\c_1\c_3-20
\c_4]\c_{n-4}+\frt14[-(21n^2-203n+472)\c_1\sp3\hs+(63n^2-609n+1430)\c_1
\c_2-(63n^2-609n+1388)\c_3]\c_{n-3}\hs+\frt14[(21n^2-203n+472)\c_1\sp2
+(63n^2-609n+1408)\c_2]\c_{n-2}\hs+\frt1{16}[105n^4-1890n^3+12131n^2
-32242n+28800]\c_1\c_{n-1}\hs+\frt1{96}n[63n^5-1575n^4+15435n^3-73801
n^2+171150n -152696]\c_n.\ea\]} \vspace{-6pt}

\n We see from \rf{K} that, as $k$ increases, $\chii^{(k)}(-1)$
involves progressively more Chern numbers, and only if $k=n$ is even
do we obtain an expression \be\chii^{(n)}(-1)=n!\chii^n= (-1)^n
n!\chii^0\ee{top} (a multiple of the Todd genus) in which the term
$\c_1\sp n$ appears. Taking $k=4$ in \rf{K},

\pro{Theorem 3} Let $M$ be a compact almost complex manifold of real
dimension $2n$ with $\c_1=0$. Then \[\left<\c_2\c_{n-2}-\c_3\c_{n-3},
[M]\right>=10\chii^{(iv)}(-1)-\frs1{24}n(15n^3-150n^2+485n-502)
\chii(-1).\]

\n For example, when $n=4$ and $\c_1=0$, this is equivalent to the fact
that the top Todd class reduces to $(3\c_2\sp2-\c_4)/720$.

\pro{Corollary} Let $M$ be a compact almost complex manifold of real
dimension $2n$ with $\c_1=0$. Then
$2n\c_n+\c_2\c_{n-2}-\c_3\c_{n-3}\equiv0$ {\ns mod 5}.

\n More generally, on any almost complex $2n$-manifold,
the indicated summands of Newton's formula \[
\underbrace{n\c_n-\s_1\c_{n-1}+\cdots-\s_{k-1}\c_{n-k+1}}+
\underbrace{\s_k\c_{n-k}+\cdots+(-1)^n\s_n}=0\] are individually zero
modulo $k+1$ if if $k+1\ge3$ is prime \cite{H53,coho}. \medbreak

Let us now turn attention to the derivation of Theorem~1 from
Theorem~2.  First suppose that $M$ is a compact \hk\ manifold of real
dimension $d=2n=4m$. If we fix a compatible complex structure $J_1$
then the closed 2-form $\eta=\w_2+i\w_3$ is holomorphic, and $\eta^m$
is nowhere zero. Wedging by $\eta^{m-p}$ is known to induce an
isomorphism $H^{p,q}\to H^{n-p,q}$ between the appropriate Dolbeault
cohomology spaces, which implies that the Hodge numbers of a
\hk\ manifold are invariant by the mirror symmetry \rf{ms}; this
was proved by Fujiki \cite{F}. Also, the $(n,0)$-form $\eta^m$
trivializes the canonical bundle of $M$, so that $\c_1=0$.

The Betti and Hodge numbers of the \ka\ manifold $M$ of real dimension
$2n$ are related by the formula $P(t)=H(t,t)$, where \[
H(x,y)=\sum_{p,q=0}^nh^{p,q}x^py^q\] is the so-called Hodge
polynomial. Because $H(x,y)$ is symmetric in $x$ and $y$, we have
\[\bar P''(-1)\a H_{xx}(-1,-1)+2H_{xy}(-1,-1) +H_{yy}(-1,-1)\\[4pt]\a
2(H_{xx}(-1,-1)+H_{xy}(-1,-1)).\ea\]  Firstly, $H(t,-1) =\chii(t)$, so
\[ H_{xx}(-1,-1)=\chii''(-1).\] Secondly, \[\bar H_{xy}(-1,-1)\a
\sum_{p,q=0}^n(-1)^{p+q}pq\hpq\\\a-\frs12n\chii'(-1)+\frs12\!
\sum_{p,q=0}^n(-1)^{p+q}[pq-(n-p)q]\hpq\\\a
\frs14n^2\chii(-1)+\frs12\!\sum_{p,q=0}^n(-1)^{p+q}pq(\hpq-h^{n-p,q}).
\ea\] The last equality uses an analogue of \rf{der} that follows from
\rf{Sd}.

Combining the above equations, we obtain

\pro{Lemma} $\ds\Phi(M)= 2\big[\,6\chii''(-1)+\frs12n(5-3n)\chii(-1)
\,\big]+6\!\sum_{p,q=0}^n(-1)^{p+q}pq(\hpq-h^{n-p,q})$.

\n From Theorem~2 we deduce that when $\c_1\c_{n-1}=0$, the symmetry
\rf{ms} reverses the sign of $\Phi$. The latter is therefore zero in
the \hk\ case.

\subsection*{3. Examples and remarks}

If $M$ and $N$ are both \hk, then their Riemannian product $M\times N$
admits an obvious \hk\ structure. We must therefore understand why the
constraint of Theorem~1 is preserved by the process of taking
products.

First suppose that $M$ is an almost complex manifold of real dimension
$2n$ with $e(M)\ne0$. If we set
\[\ga(M)=\frac{\left<\c_1\c_{n-1},[M]\right>}{\left<\c_n,
[M]\right>},\] the identity $\c(M\times N)=\c(M)\c(N)$ for the total
Chern class implies that \be\ga(M\times N)= \ga(M)+\ga(N).\ee{add} It
follows from Theorem~2 that the quantity \[
\frac{6\chii''(-1)+\frt12n(5-3n)\chii(-1)}{\chii(-1)},\]
and therefore \[\psi =
\frac{4\chii''(-1)}{\chii(-1)}-n^2,\]
also satisfies \rf{add} in place of $\ga$. When $M$ has even real
dimension $d$, the same must be true of
\be \phi = \frac{4P''(-1)}{P(-1)}-d^2\ee{phi} (also defined in
\cite{torino}), since the coefficients of $P(t)$ satisfy
\rf{Sd} with $d$ in place of $n$.

In fact the additivity of $\phi$ is as elementary as that of $\ga$,
and can deduced from the observation that $\phi$ is proportional to
the coefficient of $t^2$ in the formal power series \[\log P(-1+t) =
\log e(M) -\frs12d\,t + \frs18\phi\,t^2+\frs1{24}(3\phi+2d)t^3+
\cdots \] and similar remarks apply to $\psi$. Referring to Theorem~1,
and still assuming that $e(M)\ne0$, we have
$\Phi(M)=(3\phi(M)+5d)e(M)/2$.  Thus, $\Phi(M)$ is zero if and only if
\[ \frac{\phi(M)}{\dim(M)} = -\frs53,\] and we may think of the
right-hand side as a `coupling constant' for \hk\ manifolds.

\bigbreak\n{\bf Hilbert schemes of points.} Let $S$ be a compact complex
surface, and let $S^{(m)}$ denote its $m$-fold symmetric product
obtained by quotienting the Cartesian product by the group of
permutations.  There is a resolution \[\vep\colon S^{[m]}\lra
S^{(m)}\] in which $S^{[m]}$ is the Hilbert scheme of closed
$0$-dimensional subschemes of length $m$ on $S$, and is a smooth
complex $2m$-dimensional manifold. Each non-trivial fibre
$\vep^{-1}(x)$ is a product $(V_2)^{\alpha_2}\times\cdots\times
(V_m)^{\alpha_m}$, where $V_i=\hbox{Hilb}^i(\C[x,y])$ is the scheme
that parameterizes ideals in $\C[x,y]$ of colength $i$, and $\alpha_i$
denotes the number of $i$-tuples of points that have coallesced in
$x\in S^{(m)}$
\cite{ES,G}.

Following an example of Fujiki \cite{Fuj}, Beauville has proved
\cite{Bea} that if $S$ has a holomorphic symplectic structure then so
does $S^{[m]}$ for all $m\ge2$; its holomorphic 2-form is induced from
a natural one on $S^{(m)}$. It follows that if $S$ admits a \hk\
metric then so does $S^{[m]}$ for any $m\ge2$, and we may apply this
construction when $S$ is a K3 surface or a torus. If $S=T$ is a torus,
then $T^{[m]}$ is not locally irreducible and has a $4^m$-fold
covering of the form \rf{cov} with $k=2$ and unique non-flat factor
$Y_1$ of real dimension $4m-4$. The latter can also be viewed as a
submanifold of $T^{[m]}$ and is denoted $K_{m-1}$ in \cite{Bea} and
$A^{[\kern-.5pt[m]\kern-.5pt]}$ in \cite{GS}; when $m=2$ it is merely
the Kummer surface associated to $T$.

For any manifold $S$, the Betti numbers of $S^{(m)}$ were computed by
Macdonald \cite{Mac}. When $S$ has real dimension 4, we have
\be\sum_{m\ge0}P(S^{(m)};t)x^m=\frac{(1+tx)^{b_1}(1+t^3x)^{b_3}}
{(1-x)^{b_0}(1-t^2x)^{b_2}(1-t^4x)^{b_4}}.
\ee{Mac} The form of the right-hand side indicates its generalization
to higher dimensions, though in the present context we have $b_1=b_3$
and we may assume that $b_0=1=b_4$. Deeper results of G\"ottsche and
Soergel \cite{G,GS} give the Betti numbers of $S^{[m]}$, at least when
$S$ is a projective surface, and building on \rf{Mac}, we have \be
P(S^{[m]};t)= \sum_\alpha\>\prod_{i=1}^mP(S^{(\alpha_i)};t)
t^{2m-2\alpha_i}.\ee{GS} The sum is over all partitions $\alpha$ of
$m$, each of which is uniquely determined by $m$ non-negative integers
$\alpha_1,\ldots,\alpha_m$ with the property that $m=\sum_{i=1}^m
i\alpha_i$.

Given that the Betti numbers of $K^{[2]}$ and $K_2$ both satisfy the
constraint $\Phi=0$ given in \rf{8}, one may deduce independently of
\rf{GS} that \[\ba{c} P(K^{[2]};t)= 1+23t^2+276t^4+23t^6+t^8,\\[3pt]
 P(K_2;t) = 1+7t^2+8t^3+108t^4+8t^5+7t^6+t^8,\ea\] and the Euler
characteristics are
\[\ba{c} e(K^{[2]})=324, \\[2pt] e(K_2)=108.\ea\]

One needs \rf{GS} and a related formula in \cite{GS} to treat the
higher-dimensional cases, and it is amusing to note that \[\ba{c}
e(K^{[8]})=30178575,\\[2pt] e(K_8)=9477\ea\] are both odd. As
predicted by Theorem~1, we have
\be\Phi(K^{[m]})=0=\Phi(K_m)\ee{val} for all $m\ge1$; by contrast,
$\Phi(K^{(m)})=24m(m-1)/25$. Underlying \rf{val} is the fact that for
any projective surface $S$ with $e(S)\ne0$, one has
\be\ba{c}\phi(S^{[m]})= m \kern1pt\phi(S), \\[3pt]
\psi(S^{[m]})=m\kern1pt\psi(S).\ea\ee{realm} In particular, the
 second formula is encoded in the Hodge polynomial of $S^{[m]}$
computed in \cite{GS} and \cite{Che}; further details can be found in
\cite{coho}. We emphasize that the equations
\rf{realm} have taken us outside the realm of manifolds with $\c_1=0$,
and we shall now take the opportunity to move further afield.

\bigbreak\n{\bf Other holonomy groups.} Let $M$ be a compact oriented
manifold of real dimension $d=2n$. The constituents $P(M;-1)$ and
$P''(M;-1)$ of $\Phi(M)$ are both zero if and only if $(1+t)^4$ is a
factor of $P(M;t)$. This is certainly the case if $k\ge2$ in the
decomposition \rf{cov}. Now suppose that $M=N\times S^1$, so that \[
P(M;t)=P(N;t)(1+t),\] and the Euler characteristic $P(N;-1)$ is
zero. Then $\Phi(M)=0$ if and only if $P(N;t)$ is divisible by
$(1+t)^3$, which is equivalent to the equation $P'(N;-1)=0$.  Whilst
these observations are elementary and of little interest in the case
in which $M$ is
\hk, it seems that when $d=8$ both $\Phi(M)$ and \be
P'(N,-1)= -b_3+3b_2-5b_1+7\ee{7} are natural quantities to consider in
the context of other holonomy reductions. Now, if $M$ is an
irreducible Riemannian $d$-manifold with zero Ricci tensor whose
holonomy group $\hol(M)$ is a proper subgroup of $\SO(d)$, then
Berger's theorem \cite{Ber,Be} implies that exactly one of the
following situations occurs: (i) $M$ is \ka\ and
$\hol(M)\subseteq\SU(d/2)$, (ii) $d=7$ and $\hol(M)\cong \G_2$, or
(iii) $d=8$ and $\hol(M)\cong\spin$.

Various examples of compact 7-manifolds of type (ii) have now been
described by Joyce \cite{J} by smoothing orbifolds of the form
$T^7/\Ga$, where $\Ga$ is a finite group acting on $T^7$ preserving
the flat $\G_2$-structure. The latter is characterized by an invariant
3-form $\varphi=\sum_{i=1}^7\varphi_i$ on $\R^7$ which is the sum of
simple 3-forms $\varphi_i$, $1\le i\le7$. The first example announced
by Joyce actually had Betti numbers \be b_1=0,\quad b_2=12,\quad
b_3=43,\ee{list} so that \rf{7} vanishes. In this case,
$\Ga\cong(C_2)^3$ is a abelian group of affine transformations
preserving $\varphi$, and each element of order 2 acts trivially on
the 3-dimensional subspace of $\R^3$ determined by exactly one of the
$\varphi_i$; moreover $T^7/\Ga$ has Betti numbers $b_1=0=b_2$, $b_3=7$
so that $P'(T^7/\Ga;-1)=0$. Topologically, the smoothing process
replaces each of 12 singular 3-tori by $T^3\times S^2$; each
replacement adds 1 to the second Betti number and 3 to the third,
thereby preserving \rf{7}. Other examples turned out to be less
respectful of \rf{7}, although at the orbifold level for any $\Ga$ it
is necessarily the case that $P'(T^7/\Ga;-1)\ge0$.

If $N$ has type (ii) above, then the holonomy of the product metric on
$M=N\times S^1$ is certainly a subgroup of \spin. Now, manifolds with
holonomy contained in \spin\ generalize those with holonomy equal to
$\SU(4)$, and at the same time the two structure groups have much in
common. For example, the equation $\left<4\p_2-\p_1\sp2,[M]\right>=
8e(M)$, which is valid for any almost complex 8-manifold $M$ with
$\c_1=0$ (since then $\p_1=-2\c_2$ and $\p_2=2\c_4+\c_2\sp2$), is in
fact also satisfied in the $\spin$ case, essentially because $\SU(4)$
and $\spin$ share a maximal torus \cite{GG}. Given that $\SU(4)$
mirror symmetry changes the sign of $\Phi$, in addition to fixing $b_3$
and $2b_2+b_4$, it is curious to see what can be said about these
quantities in the \spin\ case. For example, according to Joyce, there
is a compact 8-manifold with holonomy equal to \spin\ whose Betti
numbers \be b_1=0,\quad b_2=12,\quad b_3=16,\quad b_4=150\ee{81}
satisfy the constraint $\Phi=0$ (see \rf{8}). In this case, some
standard index theory yields in addition \[ b_4^-=3b_2+7=43\] (cf.\
\rf{list}), where $b_4^-$ is the dimension of the space of harmonic
anti-self-dual 4-forms $\alpha$ (i.e., those satisfying
$*\alpha=-\alpha$).

The index theory calculations in Section~2 grew out of similar ones
for the class of quaternion-\ka\ manifolds, which are characterized by
the condition $\hol(M)\subseteq\Sp(d/4)\Sp(1)$. Unless
$\hol(M)\subseteq\Sp(d/4)\subseteq SU(d/2)$ (the \hk\ case), there is
no compatible
\ka\ structure and the Ricci tensor is non-zero.  Nevertheless the
Betti numbers of a compact quaternion-\ka\ manifold $M$ of positive
scalar curvature are also subject to a linear relation analogous to
Theorem~1. It is known that the odd Betti numbers of $M$ are zero, and
that if $d=4m$ then the integers \[
\beta_{2k}=b_{2k}(M)-b_{2k-4}(M),\qquad 0\le k\le m,\] (where
$b_j=0$ for $j<0$) are non-negative \cite{F}.  These `primitive Betti
numbers' are subject to their own remarkable constraint \be
\sum_{k=1}^m k(m+1-k)(m+1-2k)\beta_{2k}=0,\ee{ss} which is an
equivalent but neater form of the relation on the ordinary Betti
numbers of $M$ that appears in
\cite{LS,Opava}.

The character of \rf{ss} no doubt reflects its validity for certain
symmetric spaces which are the obvious (and conceivably only)
candidates for $M$. For instance, if $m=7$ then \rf{ss} becomes \[
7\beta_2+8\beta_4+5\beta_6=5\beta_{10}+ 8\beta_{12}+ 7\beta_{14},\]
and is satisfied not only by the projective space ${\Bbb H\Bbb P}^7$
($\beta_{2k}\equiv0$), the Grassmannians ${\Bbb G}{\ns r}_2(\C^9)$
($\beta_{2k}\equiv1$) and $\widetilde{{\Bbb G}{\ns r}}_4(\R^{11})$
($\beta_{4k}\equiv1$, $\beta_{4k-2}\equiv0$), but also in a striking
way by the exceptional space $\hbox{\ns F}_4/\Sp(3)\Sp(1)$
($\beta_{2k}=0$ for $k\ne4$, $\beta_8=1$). The arithmetic is more
intriguing for the $\ns E$-type quaternion-\ka\ symmetric spaces,
since their primitive Betti numbers satisfy the equation
$\beta_{2m-2k}=\beta_{2k}$ which is not obviously compatible with the
$k\leftrightarrow m+1-k$ invariance of \rf{ss}! The resulting theory
is closely related to that of the signature and other elliptic genera
of homogeneous spaces, which is the subject of
\cite{HS}.

A natural notion of self-dual connection for vector bundles over \hk\
and quaternion-\ka\ manifolds has long been known. This theory extends
in a natural way to any Riemannian manifold with reduced holonomy, by
requiring that the curvature of the connection take values in a simple
summand of the holonomy algebra. For each such algebra, the resulting
problem has an well-defined elliptic complex whose first cohomology
group parameterizes infinitesimal deformations of the connection
\cite{RRC}. Investigation of the corresponding moduli spaces can be
expected to provide further interplay between index theory and the
topology of the manifolds discussed above.

\bigbreak\n{\bf An example with $\c_1>0$.} We conclude with an example
of a \ka\ manifold with ample anti-canonical bundle. Namely, let
$\M_g$ denote the moduli space of stable rank 2 vector bundles $V$
over a Riemann surface of genus $g\ge2$ with $\ext2V$ isomorphic to a
fixed line bundle of degree one.  Then $\M_g$ is a smooth Fano
manifold of complex dimension $n=3g-3$ and index 2, its first Chern
class $\c_1$ being twice a positive integral class \cite{R}. The Betti
numbers of $\M_g$ were first given by Newstead \cite{N}, and
subsequently determined by a variety of methods \cite{AB,DR,HN}. In
particular, Atiyah and Bott include a comparison of their own
equivariant Morse theory methods with the number-theoretic approaches,
and also explain how $\M_g$ arises from an infinite-dimensional
symplectic quotient construction. The Poincar\'e polynomial of $\M_g$
is\[P(t)=\frac{(1+t^3)^{2g}-t^{2g}(1+t)^{2g}}{(1-t^2)(1-t^4)} =
(1+t)^{2g-2}\sum_{i=0}^{g-1}(1-t+t^2)^{2i} t^{2g-2-2i}.\] In
particular, for all $g$, \[ b_2(\M_g)=1, \quad b_3(\M_g)=2g, \] the
sum $P(1)$ of the Betti numbers equals $2^{2g-2}g$, and the Euler
characteristic $P(-1)$ is zero.

Because of the property $P^{(i)}(-1)=0$, $i\le2g-1$, the
simply-connected manifold $\M_g$ obviously satisfies the constraint
$\Phi(\M_g)=0$ of Theorem~1. This leads one to consider the polynomial
$\chii_t=\chii(t)$ which can in theory be computed from a knowledge of
the Chern classes $\c_k$ of $\M_g$. Newstead and Ramanan made a number
of conjectures concerning the characteristic classes of $\M_g$ which
have been subsequently proved. In particular, $\c_k=0$ for $k>2g-2$
\cite{G}, and the Pontrjagin ring (which is known to be generated
solely by $\p_1$) vanishes in degrees $4g$ and above
\cite{Ki,Th}. We shall illustrate consequences of these facts in the
relatively simple case $g=3$ which is nonetheless indicative of the
general situation.

The relations $\c_5=0=\c_6$ on $\M_3$ imply immediately that
\[\chii(-1+t)=at^4+bt^5+ct^6,\] where $a,b,c\in\Z$. From \rf{top},
$c=\chii^0$ is the Todd genus, and equals 1 because
$\c_1>0$. Furthermore, $a=\chii^{(iv)}(-1)/4!$ is given by
\[\left<\K_{6,4}^\bullet,[\M_3]
\right>=\frs1{720}\big[\,(-\c_1^3+3\c_1\c_2-3\c_3)\c_3+(\c_1^2+
3\c_2)\c_4\,\big],\] which is seen to equal 4, for example using
expressions for the Chern classes in
\cite{R}. Consequently,\[\chii(t)=(1+t)^4((b+5)+ (b+2)t+t^2),\] and
from \rf{Sd} we deduce that $b=-4$. We thereby arrive at a special
case of the more general result\be\chii(t)=(1+t)^{2g-2}
(1-t)^{g-1}\ee{chit} which was proved in \cite{NR}, and could also be
derived from a knowledge of the Hodge polynomial of $\M_g$. Given that
$b_2(\M_g)=1$, $e(\M_g)=0$, and $\c_1\ne0$, the vanishing of
$\c_{3g-2}$ may be deduced from that of $\chii''(-1)$, although the
vanishing of $\c_k$ for $2g-2<k<3g-2$ appears to be altogether
deeper. Finally, observe that \rf{chit} implies that the signature
$\chii(1)$ of $\M_g$ is zero for all $g$, which is consistent with the
vanishing of all Pontrjagin numbers.

\bigbreak\footnotesize\n{\sl Acknowledgments.} The author thanks
S.~Donaldson, F.~Hirzebruch, and M.~Thaddeus for useful
conversations. He is also grateful to R.~Jung who computed the formula
for $\K_{n,6}^\bullet$ in Section~2, and to D.~Joyce for information
provided in Section~3.

\renewcommand{\thebibliography}{\list{\arabic{enumi}.\hfil}
{\settowidth\labelwidth{18pt}\leftmargin\labelwidth\advance
\leftmargin\labelsep\usecounter{enumi}}\def\newblock{\hskip.05em}
\sloppy\sfcode`\.=1000\relax}\newcommand{\bi}{\vspace{-4pt}\bibitem}

\end{document}